\documentclass{doublecol}

\usepackage{fancyvrb}
\usepackage{psfrag}
\usepackage{algorithm}
\usepackage{graphicx}
\usepackage{amssymb}
\usepackage{moreverb,relsize}
\usepackage{amsbsy,natbib}
\usepackage{amsthm}
\usepackage{amsfonts,amsmath,bm}
\usepackage[dvipsone]{epsfig}
\usepackage{url}
\usepackage{algorithm}
\usepackage{color}

\DefineVerbatimEnvironment{code}{Verbatim}{frame=single,rulecolor=\color{blue}}

\newcommand{\tab}{\hspace*{2em}}
\newcommand{\R}{\mathbb{R}}
\newcommand{\codesize}{\footnotesize}

\begin{document}

\LRH{Unified Form-assembly Code}
\RRH{Anders Logg}

\VOL{X}
\ISSUE{X/X/X}
\PUBYEAR{2008}

\setcounter{page}{1}
\BottomCatch

\title{Efficient Representation of Computational Meshes}

\authorA{A. Logg}
\affC{Center for Biomedical Computing, Simula Research Laboratory \\
Department of Informatics, University of Oslo \\
E-mail: logg@simula.no\thanks{This research is supported by an
Outstanding Young Investigator grant from the Research Council of Norway, NFR 180450.}}

\begin{abstract}
  We present a simple yet general and efficient approach to
  representation of computational meshes. Meshes are represented as
  sets of \emph{mesh entities} of different topological dimensions and
  their \emph{incidence relations}. We discuss a straightforward and
  efficient storage scheme for such mesh representations and efficient
  algorithms for computation of arbitrary incidence relations from a
  given initial and minimal set of incidence relations. The general
  representation may harbor a wide range of computational meshes, and
  may also be specialized to provide simple user interfaces for
  particular meshes, including simplicial meshes in one, two and three
  space dimensions where the mesh entities correspond to vertices,
  edges, faces and cells.  It is elaborated on how the proposed
  concepts and data structures may be used for assembly of variational
  forms in parallel over distributed finite element meshes.
  Benchmarks are presented to demonstrate efficiency in terms of CPU
  time and memory usage.
\end{abstract}

\KEY{Mesh, mesh representation, mesh algorithms, mesh entity,
parallel assembly}

\REF{to this paper should be made as follows:
A. Logg (2008)
`Efficient Representation of Computational Meshes'.
An early version of this manuscript was presented at
MekIT´07: Fourth national conference on Computational Mechanics.
}

\maketitle

\section{Introduction}

The computational mesh is a central component of any software
framework for the (mesh-based) solution of partial differential
equations. To reduce run-time and enable the solution of large
problems, it is therefore important that the computational mesh may be
represented efficiently, both in terms of the speed of operations on
the mesh or access of mesh data, and in terms of the memory usage for
storing any given mesh in memory.

It is furthermore important that the data structure for the
representation of the mesh is general enough to harbor a wide range of
computational meshes. This generality must also be reflected in the
programming interface to the mesh representation, to allow the
implementation of general algorithms on the computational mesh. Many
algorithms, such as the assembly of a linear system from a finite
element variational problem may be implemented similarly for
simplicial, quadrilateral and hexahedral meshes if the programming
interface to the mesh representation does not enforce a specific
interface limited to a specific mesh type. For example, if the
entities on the boundary of a mesh (the \emph{facets}) may be accessed
in a similar way independently of the mesh dimension and not as
\emph{edges} in two space dimensions and \emph{faces} in three space
dimensions, one may use the same code to apply boundary conditions in
2D and 3D.

In~\cite{KnepleyKarpeev07A}, a very general and flexible
representation of computational meshes is presented. The mesh is
represented as a \emph{sieve}, which is in general a directed acyclic
graph with the mesh entities as points and directed edges describing
how the mesh entities are connected. In this paper, we take a slightly
less general approach but build on some of the concepts
from~\cite{KnepleyKarpeev07A}. In particular, we will represent the
mesh as a set of \emph{mesh entities} (corresponding to the
\emph{points} of the sieve) and their \emph{incidence relations}. We
also acknowledge the works~\cite{Ber02,Ber06}, where similar concepts
are defined and where the importance of \emph{mesh iterators} for
expressing generic algorithms on computational meshes is advocated.

The data structures and algorithms discussed in this paper have been
implemented as a C++ library and is distributed as part of DOLFIN,
see~\cite{www:dolfin}. DOLFIN is a problem-solving environment for
ordinary and partial differential equations and is developed as part
of the FEniCS project for the automation of computational mathematical
modeling, see ~\cite{www_fenics,logg:preprint:10}. Interfaces to
DOLFIN are available in the form of a C++ and a Python class library.

\subsection{Design Goals}

When designing the mesh library, we had the following design goals in
mind for the mesh representation and its interface. The mesh
representation should be \emph{simple}, meaning that the data is
represented in terms of basic C++ arrays \texttt{unsigned int*} and
\texttt{double*}; it should be \emph{generic}, meaning that it should not be
specialized to say simplicial meshes in one, two and three space
dimensions; and it should be \emph{efficient}, meaning that operations on the
mesh or access of mesh data should be fast and the storage should
require minimal memory usage for any given mesh. Furthermore, the
programming interface to the mesh representation should be
\emph{intuitive}, meaning that suitable abstractions (classes) should
be available, including specialized interfaces for specific types of
meshes as well as generic interfaces that enable dimension-independent
programming; and it should be \emph{efficient}, meaning that the
overhead of the object-oriented interface should be minimized.

\subsection{Outline}

In the following section, we present the basic concepts that define
the mesh representation and its interface. We then discuss the data
structures of the C++ implementation of the mesh representation in
DOLFIN, followed by a discussion of the algorithms used in DOLFIN to
compute any given incidence relation from a given minimal set of
incidence relations. Next, we demonstrate the programming interface to
the mesh library. This is followed by a discussion of distributed
(parallel) mesh data structures. Finally, we present a series of
benchmarks to demonstrate the efficiency of the mesh representation
and its implementation followed by some concluding remarks.

\section{Concepts}

The mesh representation is based on the following basic concepts:
\emph{mesh}, \emph{mesh topology}, \emph{mesh geometry}, \emph{mesh
entity} and \emph{mesh connectivity}. Each of these concepts is mapped
directly to the corresponding component (class) of the implementation.

\begin{figure}[htbp]
  \begin{center}
    \includegraphics[width=8.5cm]{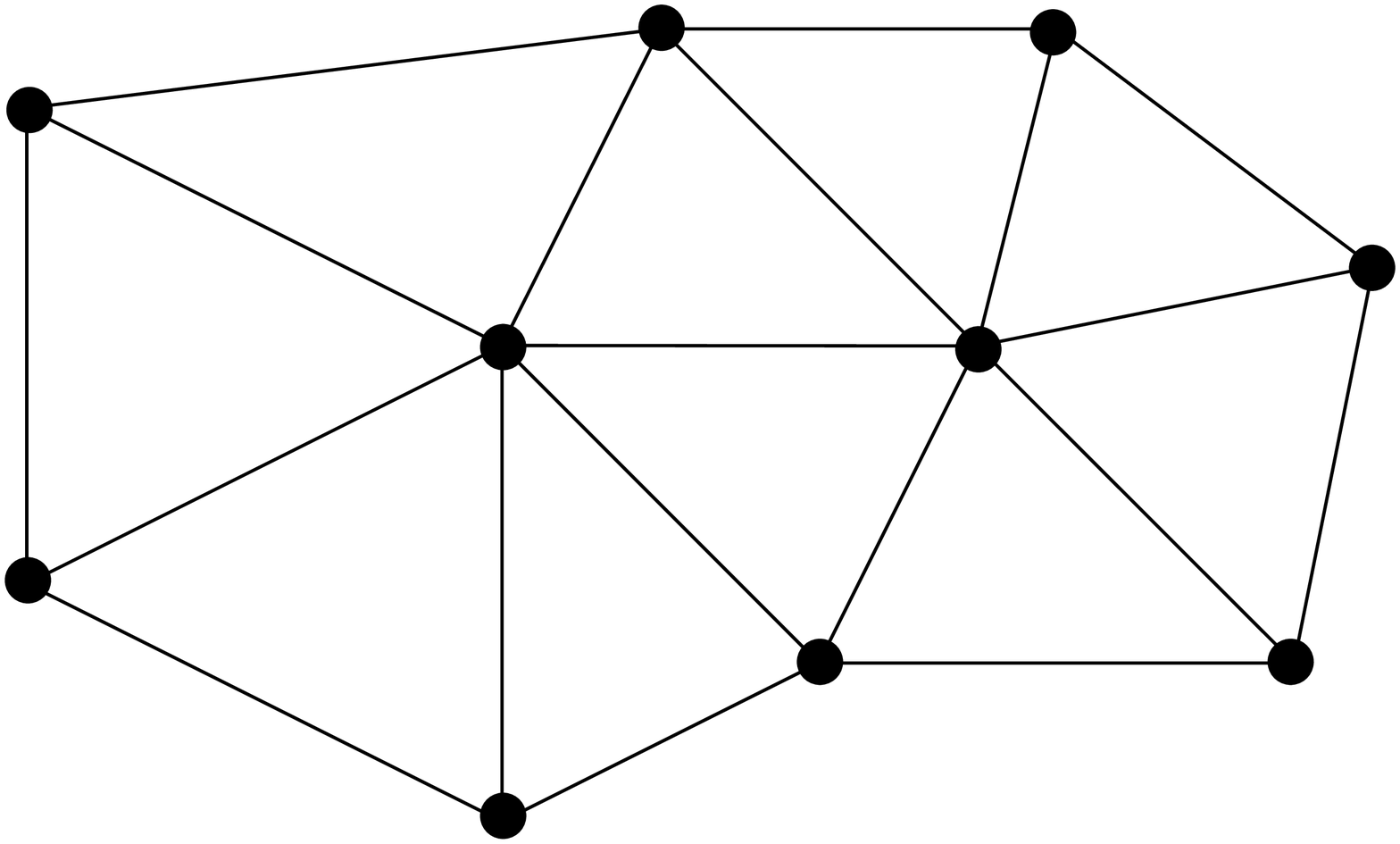}
    \caption{A mesh topology is a set of mesh entities (vertices,
    edges, etc.) and their connectivity (incidence relations), that
    is, which entities are connected (incident) to which entities.}
    \label{fig:mesh}
  \end{center}
\end{figure}

A mesh is defined by its topology and its geometry. The mesh topology
defines how the mesh is composed of its parts (the mesh entities) and
the mesh geometry describes how the mesh is embedded in some metric
space, typically $\R^n$ for $n = 1,2,3$. A mesh topology
(Figure~\ref{fig:mesh}) may be specified as a set of mesh entities
(the vertices, edges etc.) and their connectivity (incidence
relations).  Different embeddings (geometries) may be imposed on any
given mesh topology to create different meshes, e.g., when moving the
vertices of a mesh in an ALE computation. Below, we discuss the two
basic concepts mesh entity and mesh connectivity in some detail and
also introduce the concept \emph{mesh function}.

\subsection{Mesh Entities}

\begin{figure}[htbp]
  \begin{center}
    \includegraphics[width=8.5cm]{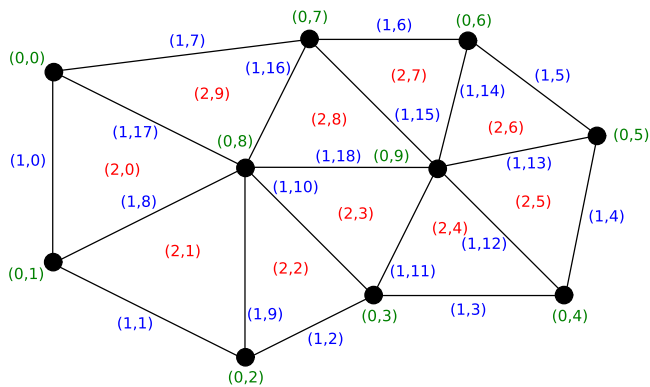}
    \caption{Each mesh entity of a mesh is identified with a pair
      $(d, i)$, where $d$ is the topological dimension of the mesh
      entity and where $i$ is a unique index for the mesh entity
      within its topological dimension, ranging from $0$ to $N_d - 1$
      with $N_d$ the number of entities of topological dimension $d$.}
    \label{fig:meshentities}
  \end{center}
\end{figure}

A mesh entity is a pair $(d, i)$, where $d$ is the topological
dimension of the mesh entity and where $i$ is a unique index for the
mesh entity within its topological dimension, ranging from $0$ to $N_d
- 1$ with $N_d$ the number of entities of topological dimension
$d$. We let $D$ denote the maximal topological dimension over the mesh
entities and set the topological dimension of the mesh equal to $D$.
This is illustrated in Figure~\ref{fig:meshentities}, where each mesh
entity is labeled by its topological dimension and index $(d, i)$.

For convenience, we also name common entities of low topological
dimension or codimension. We refer to entities of topological
dimension~0 as \emph{vertices}, entities of dimension~1 as
\emph{edges}, entities of dimension~2 as \emph{faces}, entities of
codimension~$1$ (dimension $D-1$) as \emph{facets} and entities of
codimension~$0$ (dimension $D$) as \emph{cells}. Thus, for a
triangular mesh, the edges are also facets and the faces are also
cells, and for a tetrahedral mesh, the faces are also facets. This is
summarized in Table~\ref{tab:namedentities}.

\begin{table}[htbp]
  \begin{center}
    \begin{tabular}{|l|c|c|}
      \hline
      Entity & Dimension & Codimension \\
      \hline
      \hline
      Vertex     & $0$ &     $D$ \\
      \hline
      Edge       & $1$ &     $D - 1$ \\
      \hline
      Face       & $2$ &     $D - 2$ \\
      \hline
      Facet      & $D - 1$ & $1$ \\
      \hline
      Cell       & $D$ &     $0$ \\
      \hline
    \end{tabular}
    \caption{Named entities of low topological dimension or codimension.}
    \label{tab:namedentities}
  \end{center}
\end{table}

\subsection{Mesh Connectivity}

We refer to the set of incidence relations on a set of mesh entities
as the \emph{connectivity} of the mesh. For a mesh of topological dimension
$D$, there are $(D+1)^2$ different classes of incidence relations
(connectivities) to consider. Each such class is denoted here by $d
\rightarrow d'$ for $0 \leq d, d' \leq D$. For any given mesh entity
$(d, i)$, its connectivity $(d \rightarrow d')_i$ is given by the set of
incident mesh entities of dimension $d'$.

Thus, for a triangular mesh (of topological dimension $D = 2$), there
are nine different incidence relations of interest between the
entities of the mesh. These are in turn $0 \rightarrow 0$ (the
vertices incident to each vertex), $0 \rightarrow 1$ (the
edges incident to each vertex), \ldots, $D \rightarrow D$ (the
cells incident to each cell).

For $d > d'$, the definition of incidence is evident. Mesh entity
$(d', i')$ is incident to mesh entity $(d, i)$ if $(d', i')$ is
\emph{contained} in $(d, i)$. Thus, the three vertices of a triangular
cell form the set of incident vertices and the three edges form the
set of incident edges. For $d < d'$, we define mesh entity $(d', i')$
as incident to mesh entity $(d, i)$ if $(d, i)$ is incident to $(d',
i')$. It thus remains to define incidence for $d = d'$. For $d, d' >
0$, we say that mesh entity $(d', i')$ is incident to mesh entity $(d,
i)$ if both are incident to a common vertex, that is, a mesh entity of
dimension zero, while for $d = d' = 0$, we say that $(d', i')$ is
incident to $(d, i)$ if both are incident to a common cell, that is, a
mesh entity of dimension $D$.

Together, the set of mesh entities and the connectivity (incidence
relations) define the topology of the mesh. Note that the complete set
of incidence relations $d \rightarrow d'$ for $0 \leq d, d' \leq D$
may be determined from the single class of incidence relations $D
\rightarrow 0$, that is, the vertices of each cell in the mesh. We
return to this below when we present an algorithm for computing any
given class of incidence relations from the minimal set of incidence
relations $D \rightarrow 0$.

\subsection{Mesh Functions}

We define a \emph{mesh function} as a discrete function that takes a
value on the set of mesh entities of a given fixed dimension $0 \leq d
\leq D$. Mesh functions are simple objects but very useful. A
real-valued mesh function may for example be used to describe material
parameters on the cells of a mesh. A boolean-valued mesh function may
be used to set markers on cells or edges for adaptive
refinement. Integer-valued mesh functions may be used to express
inter-connectivity between two separate meshes. A typical use is when
a boundary mesh is extracted from a given mesh (by identifying the set
of facets that are incident to exactly one cell). One may then use a
mesh function to describe the mapping from the cells in the extracted
boundary mesh (which has topological dimension $D - 1$) to the
corresponding facets in the original mesh (which has topological
dimension $D$). Note that mesh functions are discrete and are not
meant to represent for example a piecewise polynomial finite element
function on the mesh.

\section{Data Structures}

The mesh representation as described in the previous section has been
implemented as a small C++ class library and is available freely as
part of the DOLFIN C++ finite element library, version 0.6.3 or
higher. Each of the basic concepts \emph{mesh}, \emph{mesh topology},
\emph{mesh geometry}, \emph{mesh entity}, \emph{mesh connectivity} and
\emph{mesh function} is realized by the corresponding class
\texttt{Mesh}, \texttt{MeshTopology}, \texttt{MeshGeometry},
\texttt{MeshEntity}, \texttt{MeshConnectivity} and
\texttt{MeshFunction}. All basic data structures are stored as static
arrays of unsigned integers (\texttt{unsigned int*}) or floating point
values (\texttt{double*}), which minimizes the cost of storing the
mesh data and allows for quick access of mesh data.  We discuss each
of these classes/data structures in detail below.

\subsection{The Class \texttt{Mesh}}

The class \texttt{Mesh} stores a \texttt{MeshTopology} and a
\texttt{MeshGeometry} that together define the mesh. The
\texttt{MeshTopology} and \texttt{MeshGeometry} are independent of
each other and of the \texttt{Mesh}. Although it is possible to work
with the \texttt{MeshTopology} and \texttt{MeshGeometry} separately,
they are most conveniently accessed through a \texttt{Mesh} class that
holds a pair of a matching topology and geometry.

\subsection{The Class \texttt{MeshTopology}}

The class \texttt{MeshTopology} stores the topology of a mesh as a
set of mesh entities and connectivities. For each pair of topological
dimensions $(d, d'), 0 \leq d, d' \leq D$, the class
\texttt{MeshTopology} stores a \texttt{MeshConnectivity} representing
the set of incidence relations $d \rightarrow d'$. The mesh entities
themselves need not be stored explicitly; they are stored implicitly
for each topological dimension $d$ as the set of pairs $(d, i)$ for $0
\leq i < N_d$, where $N_d$ is the number of mesh entities of
topological dimension $d$. Thus, for each topological dimension, the
class \texttt{MeshTopology} stores an (unsigned) integer $N_d$, from
which the set of mesh entities $\{(d, 0), (d, 1), \ldots, (d, N_d -
1)\}$ may be generated.

\subsection{The Class \texttt{MeshGeometry}}

The class \texttt{MeshGeometry} stores the geometry of a
mesh. Currently, only the simplest possible representation has been
implemented, where only the coordinates of each vertex are stored.
These coordinates are stored in a contiguous array
\texttt{coordinates} of size $n N_0$, where $n$ is the geometric
dimension and $N_0$ is the number of vertices.

\subsection{The Class \texttt{MeshEntity}}

The class \texttt{MeshEntity} provides a \emph{view} of a given mesh
entity $(d, i)$. The mesh entities themselves are not stored, but a
\texttt{MeshEntity} may be generated from a given pair $(d, i)$. The
class \texttt{MeshEntity} provides a convenient interface for
accessing mesh data, in particular in combination with the concept of
mesh iterators, as will be discussed in more detail below. Thus, one
may for any given \texttt{MeshEntity} access its topological
dimension~$d$, its index~$i$ and its set of incidence relations
(connected mesh entities) of any given topological
dimension~$d'$. Specialized interfaces are provided for the named mesh
entities of Table~\ref{tab:namedentities} in the form of the following
sub classes of \texttt{MeshEntity}: \texttt{Vertex}, \texttt{Edge},
\texttt{Face}, \texttt{Facet} and \texttt{Cell}.

\subsection{The Class \texttt{MeshConnectivity}}

The class \texttt{MeshConnectivity} stores the set of incidence
relations~$d \rightarrow d'$ for a fixed pair of topological
dimensions~$(d, d')$. The set of incidence relations is stored as a
contiguous \texttt{unsigned int} array \texttt{indices} of entity
indices for dimension $d'$ entities, together with an auxiliary
\texttt{unsigned int} array \texttt{offsets} that specifies the offset
into the first array for each entity of dimension $d$.\footnote{The
storage is similar to the standard compressed row storage (CRS) format
for sparse matrices, except that only the column indices need to be
stored, not the values. Also note that the two arrays
\texttt{indices} and \texttt{offsets} are private data structures of
the class \texttt{MeshConnectivity}. The user is presented with a more
intuitive interface, as will be demonstrated below.}  The size of the
first array \texttt{indices} is equal to the total number of incident
entities of dimension~$d'$ and the size of the second array
\texttt{offsets} is equal to the total number of entities of dimension
$d$ plus one.

\begin{figure}[htbp]
  \begin{center}
    \includegraphics[width=8.5cm]{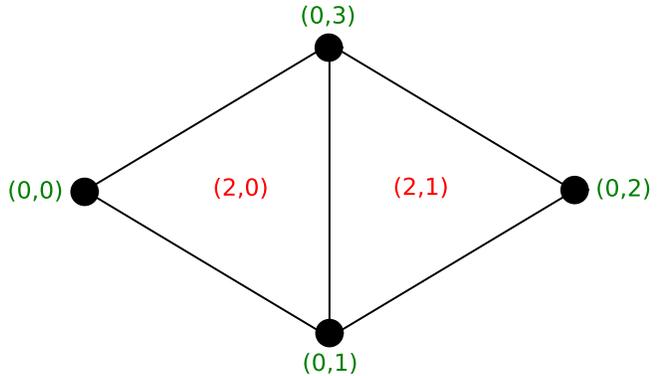}
    \caption{The mesh connectivity $2 \rightarrow 0$ (the vertices of
    each cell) for this triangular mesh with two cells and four vertices
    is stored as two arrays \texttt{indices = [0, 1, 3, 1, 2, 3]} and \texttt{offsets = [0, 3, 6]}.}
    \label{fig:connectivityexample}
  \end{center}
\end{figure}

As an example, consider the storage of the set of incidence relations
$2 \rightarrow 0$, that is the vertices of each cell, for the
triangular mesh in Figure~\ref{fig:connectivityexample}. The mesh has
two entities of dimension $d = 2$ and four entities of dimension $d' =
0$. Furthermore, each entity of dimension $d = 2$ is incident to three
entities of dimension $d' = 0$. The array \texttt{entities} is then
given by \texttt{[0, 1, 3, 1, 2, 3]} and the array \texttt{offsets} is given
by \texttt{[0, 3, 6]}.

\subsection{The Class \texttt{MeshFunction}}

The class \texttt{MeshFunction} stores a single array of $N_d$ values on the
mesh entities of a given fixed dimension $d$, and is templated over
the value type. Typical uses include \texttt{MeshFunction<double>} for
material parameters that take a constant value on each cell of a mesh,
\texttt{MeshFunction<bool>} for cell markers that indicate cells that
should be refined, and \texttt{MeshFunction<unsigned int>} to store
inter-mesh connectivity or sub domain markers.

\subsection{Minimal Storage}

The mesh data structures described above are summarized in
Table~\ref{tab:datastructures}. We note that the classes
\texttt{Mesh} and \texttt{MeshTopology} function as ``aggregate classes''
that collect mesh data stored elsewhere, and that no data is stored in
the class \texttt{MeshEntity}. All data is thus stored in the class
\texttt{MeshConnectivity} (in the two arrays \texttt{indices} and
\texttt{offsets}) and in the class \texttt{MeshGeometry} (in the array
\texttt{coordinates}). Note that one \texttt{MeshConnectivity} object
is stored for each pair of topological dimensions $(d, d')$ for which
the mesh connectivity has been initialized.

\begin{table}[htbp]
  \begin{center}
    \begin{tabular}{|l|l|}
      \hline
      Data structure        & Principal data \\
      \hline
      \hline
      \texttt{Mesh}         & \texttt{MeshTopology topology} \\
                            & \texttt{MeshGeometry geometry}  \\
      \hline
      \texttt{MeshTopology} & \texttt{MeshConnectivity**} \\
      & \hspace{1.5cm} \texttt{connectivities} \\
      \hline
      \texttt{MeshGeometry} & \texttt{double* coordinates} \\
      \hline
      \texttt{MeshEntity}   & -- \\
      \hline
      \texttt{MeshConnectivity} & \texttt{unsigned int* indices} \\
                                & \texttt{unsigned int* offsets} \\
      \hline
    \end{tabular}

    \vspace{0.5cm}

    \begin{tabular}{|l|l|}
      \hline
      Principal data        & Size \\
      \hline
      \hline
      \texttt{MeshTopology topology} & -- \\
      \texttt{MeshGeometry geometry} & -- \\
      \hline
      \texttt{MeshConnectivity** connectivities} & -- \\
      \hline
      \texttt{double* coordinates} & $nN_0$ \\
      \hline
      \texttt{unsigned int* indices} & $\mathcal{O}(N_{d})$ \\
      \texttt{unsigned int* offsets} & $N_d + 1$ \\
      \hline
    \end{tabular}
    \caption{Summary of mesh data structures.}
    \label{tab:datastructures}
  \end{center}
\end{table}

As an illustration, consider the storage of a tetrahedral mesh with
$N_0$ vertices and $N_3$ cells (tetrahedra) embedded in $\R^3$ where
we only store the set of incidence relations $D \rightarrow 0$. Each
cell has four vertices, so the class \texttt{MeshConnectivity} stores
$4N_3 + N_3 + 1 \sim 5N_3$ values of type \texttt{unsigned int}.
Furthermore, the class \texttt{MeshGeometry} stores $3N_0$ values of
type \texttt{double}. Thus, if an \texttt{unsigned int} is four bytes
and a \texttt{double} is eight bytes, then the total size of the mesh
is $20N_3 + 24N_0$ bytes. For a standard uniform tetrahedral mesh of
the unit square, the number of cells is approximately six times the
number of vertices, so the total size of the mesh is
\begin{equation}
  (20N_3 + 24N_0) \, \textrm{b} = (20N_3 + 24N_3/6) \, \textrm{b} = 24N_3 \, \textrm{b}.
\end{equation}
Thus, a mesh with $1,000,000$ cells may be stored in just $24 \, \textrm{Mb}$.
Note that if additional mesh connectivity is computed,
like the edges or facets of the tetrahedra, more memory will be
required to store the mesh.

\section{Algorithms}

In this section, we present the algorithms used by the DOLFIN mesh
library to compute the mesh connectivity $d \rightarrow d'$ for any
given $0 \leq d, d' \leq D$. We assume that we are given an initial
set of incidence relations $D \rightarrow 0$, that is, we know the
vertices of each cell in the mesh.

\begin{figure}[htbp]
  \begin{center}
    \includegraphics[width=7cm]{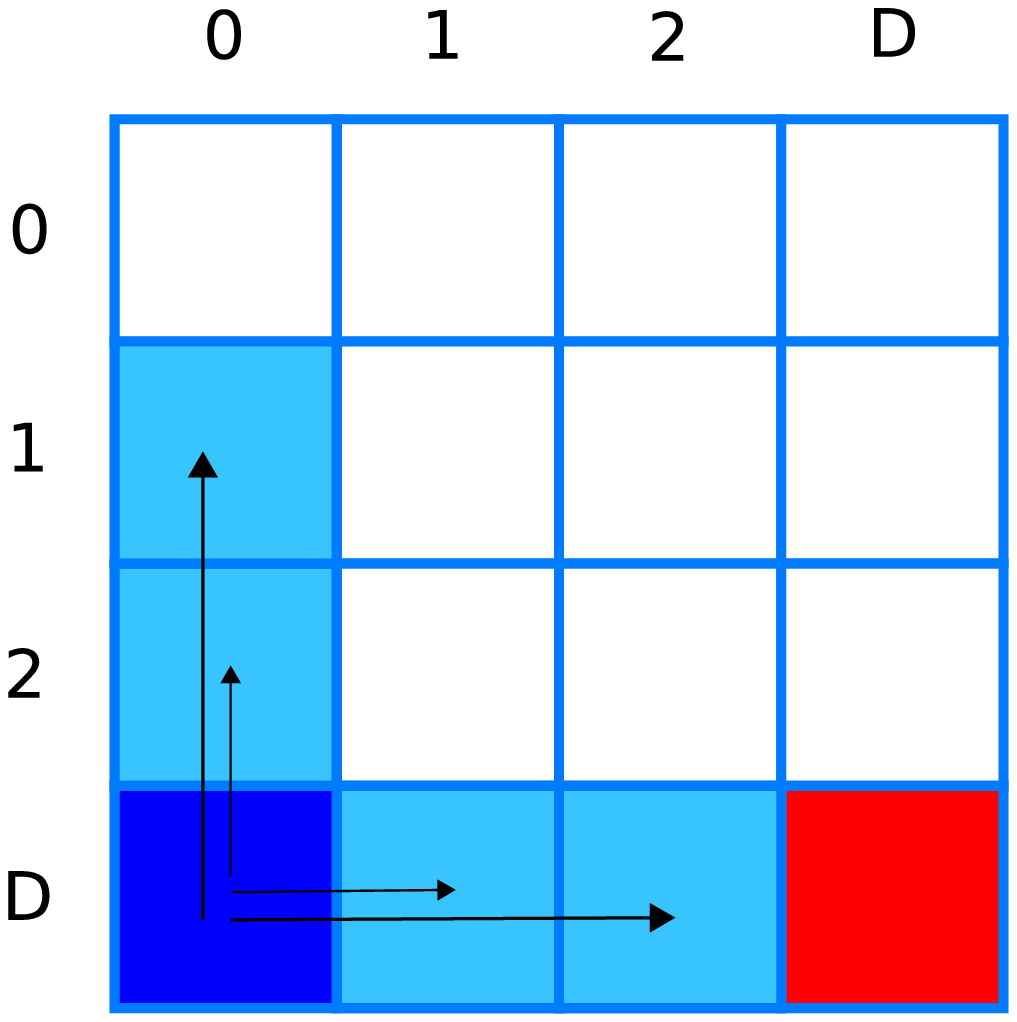} \\
    \vspace{0.5cm}
    \includegraphics[width=7cm]{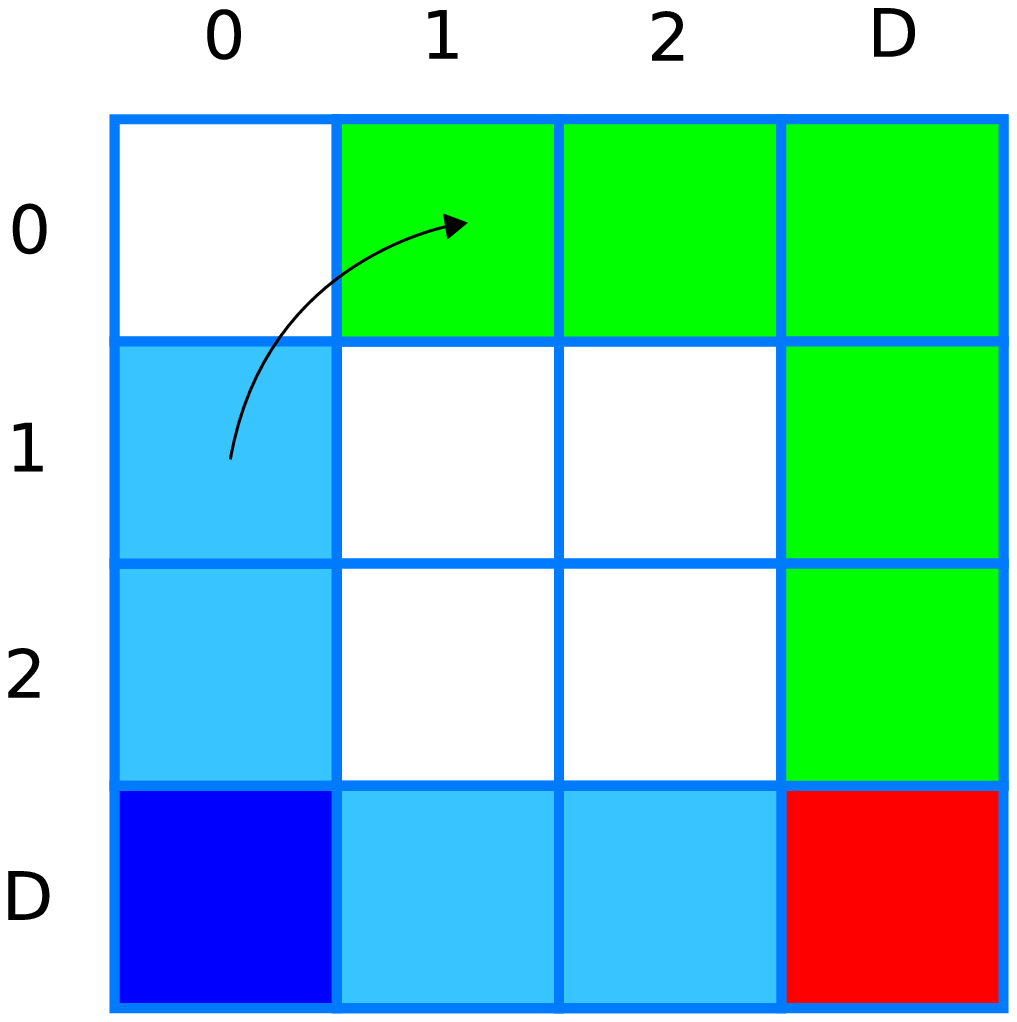} \\
    \vspace{0.5cm}
    \includegraphics[width=7cm]{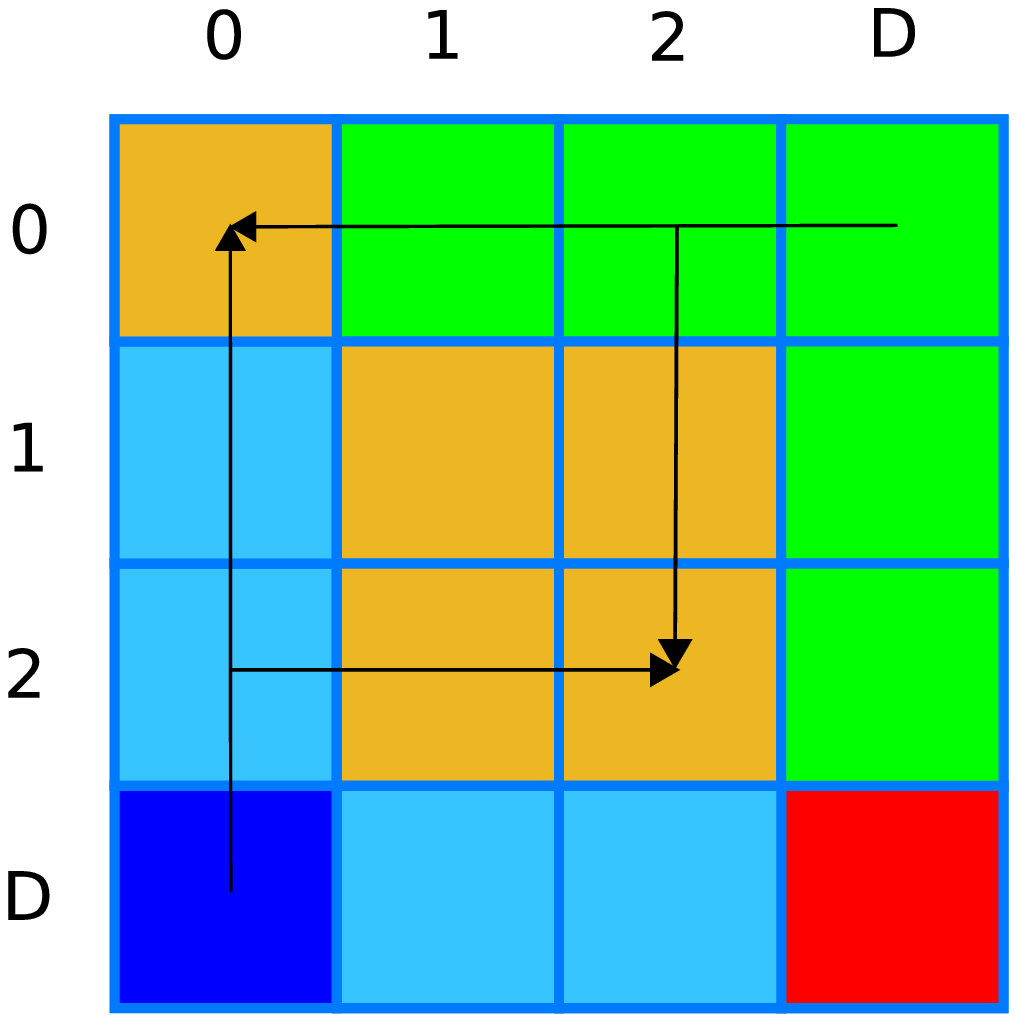}
    \caption{The three basic algorithms for computing connectivity.
      From the top:
      Build
      (computing connectivity $D \rightarrow d$ and $d \rightarrow 0$ from $D \rightarrow 0$ and $D \rightarrow D$),
      Transpose
      (computing connectivity $d \rightarrow d'$ from $d' \rightarrow d$) and
      Intersection
      (computing connectivity $d \rightarrow d'$ from $d \rightarrow d''$ and $d'' \rightarrow d'$).}
    \label{fig:algorithms}
  \end{center}
\end{figure}

The key to computing the mesh connectivities of a mesh is to compute
the connectivities in a particular order. For example, if the vertices
are known for each edge in the mesh ($1 \rightarrow 0$), then it is
straightforward to compute the edges incident to each vertex ($0
\rightarrow 1$) as will be explained below. The computation is based
on three algorithms that are used successively in a particular order
to compute the desired connectivity. As a consequence, the computation
of a certain connectivity $d \rightarrow d'$ might require the
computation of one or more other connectivities. We describe these
algorithms in detail below. An overview is given in
Figure~\ref{fig:algorithms}

\subsection{Build}

Algorithm~\ref{alg:build} (Build) computes the connectivities $D
\rightarrow d$ and $d \rightarrow 0$ from $D \rightarrow 0$ and $D
\rightarrow D$ for $0 < d < D$. In other words, given the vertices and
incident cells of each cell in the mesh, Algorithm~\ref{alg:build}
computes the entities of dimension~$d$ of each cell and for each such
entity the vertices of that entity. Thus, if $d = 1$, then the edges
of each cell and the vertices of each edge are computed.

The notation of Algorithm~\ref{alg:build} requires some
explanation. As before, we let $(d \rightarrow d')_i$ denote the
set of entities of dimension $d'$ incident to entity $(d, i)$:
\begin{equation}
  (d \rightarrow d')_i = \{(d', j) : (d', j) \mbox{ incident to } (d, i)\}.
\end{equation}
Algorithm~\ref{alg:build} also uses the operation
\begin{equation}
  d \xrightarrow{\mathrm{local} (D, i)} 0,
\end{equation}
which denotes the set of vertex sets incident to the mesh entities of
topological dimension $d$ of a given cell $(D, i)$. To make this
concrete, consider a triangular mesh (for which $D = 2$) and take $d =
1$. If $V_i = d \xrightarrow{\mathrm{local} (D, i)} 0$, then $V_i$
denotes the set of vertex sets incident to the edges of triangle
number $i$. The set $V_i$ consists of three sets of vertices (one for
each edge) and each set $v_i \in V_i$ contains two vertices.  In
addition, Algorithm~\ref{alg:build} uses the operation
\begin{equation}
  \mathrm{index}((D, j), d, v_i),
\end{equation}
which denotes the index of the entity of dimension $d$ in the cell $(D,
j)$ which is incident to the vertices $v_i$.

We may now summarize Algorithm~\ref{alg:build} as follows. For each
cell $(D, i)$, we create a set of candidate entities of dimension
$d$, represented by their incident vertices in the set $V_i$. This
operation is local on each cell and must be performed differently for
each different type of mesh. We then iterate over each cell incident
to the cell $(D, i)$ and check for each candidate entity $v_i \in V_i$
if it has already been created by any of the previously visited cells,
making sure that two incident cells agree on the index of any common
incident entity.

\begin{algorithm}[htbp]
  \begin{tabbing}
    $k = 0$ \\
    \textbf{for each} $(D, i)$ \\
    \tab $V_i = d \xrightarrow{\mathrm{local} (D, i)} 0$ \\
    \tab \textbf{for each} $(D, j) \in (D \rightarrow D)_i$ \textbf{such that} $j < i$ \\
    \tab \tab $V_j = d \xrightarrow{\mathrm{local} (D, j)} 0$ \\
    \tab \tab \textbf{for each} $v_i \in V_i$ \\
    \tab \tab \tab \textbf{if} $v_i \in V_j$ \\
    \tab \tab \tab \tab $l = \mathrm{index}((D, j), d, v_i)$ \\
    \tab \tab \tab \tab $(D \rightarrow d)_i = (D \rightarrow d)_i \cup (d, l)$ \\
    \tab \tab \tab \textbf{else} \\
    \tab \tab \tab \tab $(D \rightarrow d)_i = (D \rightarrow d)_i \cup (d, k)$ \\
    \tab \tab \tab \tab $(d \rightarrow 0)_k = v_i$ \\
    \tab \tab \tab \tab $k = k + 1$
  \end{tabbing}
  \caption{Build($d$),
    computing $D \rightarrow d$ and $d \rightarrow 0$ from $D \rightarrow 0$ and $D \rightarrow D$ for $0 < d < D$}
  \label{alg:build}
\end{algorithm}

\subsection{Transpose}

Algorithm~\ref{alg:transpose} (Transpose) computes the connectivity $d
\rightarrow d'$ from the connectivity $d' \rightarrow d$ for $d <
d'$. For each entity of dimension $d'$, we iterate over the incident
entities of dimension $d$ and add the entities of dimension $d'$ as
incident entities to the entities of dimension $d$. We may thus
compute for example the incident cells of each vertex (the cells to
which the vertex belongs) by iterating over the cells of the mesh
and for each cell over its incident vertices.

\begin{algorithm}[htbp]
  \begin{tabbing}
    \textbf{for each} $(d', j)$ \\
    \tab \textbf{for each} $(d, i) \in (d' \rightarrow d)_j$ \\
    \tab \tab $(d \rightarrow d')_i = (d \rightarrow d')_i \cup (d', j)$
  \end{tabbing}
  \caption{Transpose($d$, $d'$),
    computing $d \rightarrow d'$ from $d' \rightarrow d$ for $d < d'$}
  \label{alg:transpose}
\end{algorithm}

\subsection{Intersection}

Algorithm~\ref{alg:intersection} (Intersection) computes the
connectivity $d \rightarrow d'$ from $d \rightarrow d''$ and $d''
\rightarrow d'$ for $d \geq d'$. For each entity $(d, i)$ of dimension
$d$, we iterate over each incident entity $(d'', k)$ of dimension
$d''$ and for each such entity we iterate over each incident entity
$(d', j)$ of dimension $d'$. We then check if either $(d, i)$ and
$(d', j)$ are entities of the same topological dimension or if $(d',
j)$ is completely contained in $(d, i)$ by checking that each vertex
incident to $(d', j)$ is also incident to $(d, i)$, in which case
$(d', j')$ is added as an incident entity of entity $(d, i)$.

Here, $d''$ must be chosen according to the definition of incidence
given above. For example, we may take $d'' = 0$ to compute the connectivity $D
\rightarrow D$ (the incident cells of each cell) by iterating over the
vertices of each cell and for each such vertex iterate over the
incident cells.

\begin{algorithm}[htbp]
  \begin{tabbing}
    \textbf{for each} $(d, i)$ \\
    \tab \textbf{for each} $(d'', k) \in (d \rightarrow d'')_i$ \\
    \tab \tab \textbf{for each} $(d', j) \in (d'' \rightarrow d')_k$ \\
    \tab \tab \tab \textbf{if}
    ($d = d'$ \textbf{and} $i \neq j$) \textbf{or} \\
    \hspace{2.4cm} ($d > d'$ \textbf{and} $(d' \rightarrow 0)_j \subseteq (d \rightarrow 0)_i$) \\
    \tab \tab \tab \tab $(d \rightarrow d')_i = (d \rightarrow d')_i \cup (d', j)$
  \end{tabbing}
  \caption{Intersection($d$, $d'$, $d''$),
    computing $d \rightarrow d'$ from $d \rightarrow d''$ and $d''
    \rightarrow d'$ for $d \geq d'$}
  \label{alg:intersection}
\end{algorithm}

\subsection{Successive Application of Build, Transpose and Intersection}

Any given connectivity $d \rightarrow d'$ for $0 \leq d, d' \leq D$
may be computed by a successive application of
Algorithms~\ref{alg:build}--\ref{alg:intersection} in a suitable
order. In Algorithm~\ref{alg:connectivity}, we present the basic logic
for a successive and recursive application of the three basic
algorithms Build, Transpose and Intersection to compute any given
connectivity.

\begin{figure}[htbp]
  \begin{center}
    \includegraphics[width=8.5cm]{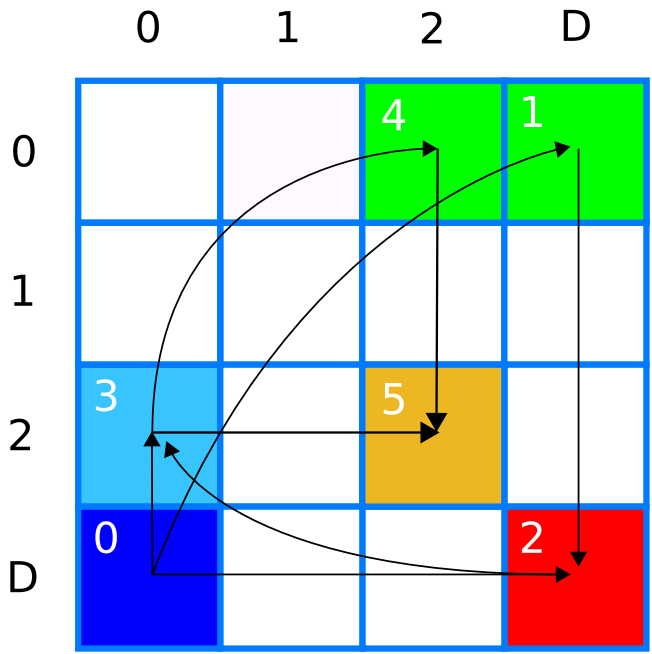}
    \caption{Computing connectivity $2 \rightarrow 2$ (the faces
    incident to any given face) by successive application of
    Transpose, Intersection, Build, Transpose and Intersection.}
    \label{fig:connectivity,all}
  \end{center}
\end{figure}

\begin{algorithm}[htbp]
  \begin{tabbing}
    \textbf{if} $N_d = 0$ \\
    \tab Build($d$) \\
    \textbf{if} $N_{d'} = 0$ \\
    \tab Build($d'$) \\
    \textbf{if} $d \rightarrow d' \neq \emptyset$ \\
    \tab \textbf{return} \\
    \\
    \textbf{if} $d < d'$ \\
    \tab Connectivity($d'$, $d$) \\
    \tab Transpose($d$, $d'$) \\
    \textbf{else} \\
    \tab \textbf{if} $d = 0$ \textbf{and} $d' = 0$ \\
    \tab \tab $d''$ = $D$ \\
    \tab \textbf{else} \\
    \tab \tab $d''$ = $0$ \\
    \tab Connectivity($d$, $d''$) \\
    \tab Connectivity($d''$, $d'$) \\
    \tab Intersection($d$, $d'$, $d''$)
  \end{tabbing}
  \caption{Connectivity($d$, $d'$),
    computing $d \rightarrow d'$ by application of Algorithms~\ref{alg:build}--\ref{alg:intersection}}
  \label{alg:connectivity}
\end{algorithm}

We illustrate this in Figure~\ref{fig:connectivity,all} for
computation of the connectivity $2 \rightarrow 2$, the incident faces
of each face, for a tetrahedral mesh. From the given connectivity $D
\rightarrow 0$, we first compute the connectivity $0 \rightarrow D$ by
an application of Transpose. This allows us to compute $D \rightarrow
D$ by an application of Intersection. The connectivity $2 \rightarrow
0$ (and $D \rightarrow 2$) may then be computed by an application of
Build. We then apply Transpose to compute $0 \rightarrow 2$ and
finally Intersection to compute $2 \rightarrow 2$.

\subsection{Memory Handling}

For each of Algorithms~\ref{alg:build}--\ref{alg:intersection}, memory
usage may be conserved by running each algorithm twice; first one
round to count the number of incident entities, which allows the static
data structures discussed above to be preallocated, and then another
round to set the values of the incident entities. Furthermore, memory
usage may be conserved by clearing incidence relations that get
computed as byproducts of
Algorithms~\ref{alg:build}--\ref{alg:intersection} when they are no
longer needed.

\section{Interfaces}

In this section, we briefly describe the user interface of the DOLFIN
mesh library. We only describe the C++ interface, but note that an
(almost) identical Python interface is also available.

\subsection{Creating a Mesh}

A mesh may be created in one of three ways, as illustrated in
Figure~\ref{fig:interface:create}. Either, the mesh is defined by a
data file in the DOLFIN XML format\footnote{A conversion script
\texttt{dolfin-convert} is provided for conversion from other popular
mesh formats (including Gmsh and Medit) to DOLFIN XML format.}, or the
mesh is defined vertex by vertex and cell by cell using the DOLFIN
mesh editor, or the mesh is defined as one of the DOLFIN built-in
meshes. Currently provided built-in meshes include triangular meshes
of the unit square and tetrahedral meshes of the unit cube.

\begin{figure}[htbp]
  \codesize
  \begin{center}
    \begin{code}
// Read mesh from file
Mesh mesh0("mesh.xml");

// Build mesh using the mesh editor
Mesh mesh1;
MeshEditor editor;
editor.open(mesh1, "triangle", 2, 2);
editor.initVertices(4);
editor.addVertex(0, 0.0, 0.0);
editor.addVertex(1, 1.0, 0.0);
editor.addVertex(2, 1.0, 1.0);
editor.addVertex(3, 0.0, 1.0);
editor.initCells(2);
editor.addCell(0, 0, 1, 2);
editor.addCell(1, 0, 2, 3);
editor.close();

// Create simple mesh of the unit cube
UnitCube mesh2(16, 16, 16);
    \end{code}
    \caption{A DOLFIN mesh may be defined either by an XML data file,
      or explicitly using the DOLFIN mesh editor, or as a built-in
      predefined mesh. The last two arguments in the call to
      \texttt{MeshEditor::open()} specify the topological and geometric
      dimensions of the mesh respectively.}
    \label{fig:interface:create}
  \end{center}
\end{figure}

\subsection{Mesh Iterators}

Mesh data may be accessed directly from the mesh, but is most
conveniently accessed through the mesh iterator interface. Algorithms
operating on a mesh (including
Algorithms~\ref{alg:build}--\ref{alg:intersection}) may often be
expressed in terms of \emph{iterators}. Mesh iterators can be used to
iterate either over the global set of mesh entities of a given
topological dimension, or over the locally incident entities of any
given mesh entity. Two alternative interfaces are provided; the
general interface \texttt{MeshEntityIterator} for iteration over
entities of some given topological dimension~$d$, and the specialized
mesh iterators \texttt{VertexIterator}, \texttt{EdgeIterator},
\texttt{FaceIterator}, \texttt{FacetIterator} and
\texttt{CellIterator} for iteration over named entities.  Iteration
over mesh entities may be nested at arbitrary depth and the
connectivity (incidence relations) required for any given iteration is
automatically computed (at the first occurrence) by the algorithms
presented in the previous section.

A \texttt{MeshEntityIterator} (\texttt{it}) may be dereferenced
(\texttt{*it}) to create a \texttt{MeshEntity}, and any member function
\texttt{MeshEntity::foo()} may be accessed by \texttt{it->foo()}. A
\texttt{MeshEntityIterator} may thus be thought of as a \emph{pointer}
to a \texttt{MeshEntity}. Similarly, the named mesh entity iterators
may be dereferenced to create the corresponding named mesh
entities. Thus, dereferencing a \texttt{VertexIterator} gives a
\texttt{Vertex} which provides an interface to access vertex data. For
example, if \texttt{it} is a \texttt{VertexIterator}, then
\texttt{it->point()} returns the coordinates of the vertex.

\begin{figure}[htbp]
  \begin{center}
    \includegraphics[width=8.5cm]{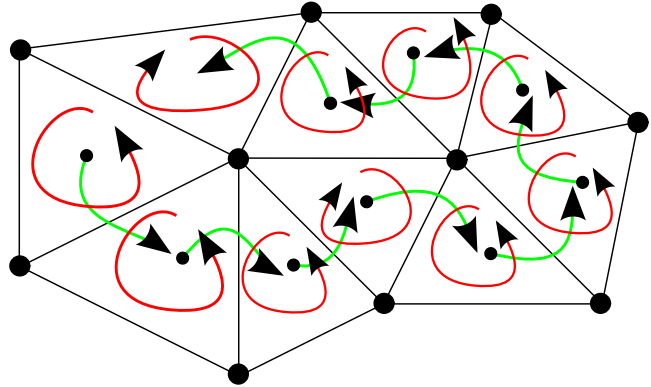}
    \caption{Iteration over all vertices of all cells in a mesh. The
    order of iteration is decided by the definition of the mesh, or
    alternatively, the UFC ordering convention~\cite{logg:www:08} if
    the mesh is ordered. Meshes may be ordered by a call to
    \texttt{Mesh::order()}.}
    \label{fig:iteration}
  \end{center}
\end{figure}

The use of mesh iterators is demonstrated in
Figure~\ref{fig:interface:iterators} for iteration over all cells in
the mesh and for each cell all its vertices as illustrated in
Figure~\ref{fig:iteration}. For each cell and each vertex, we print
its mesh entity index. We also demonstrate the use of named mesh
entity iterators to print the coordinates of each vertex.

\begin{figure}[htbp]
  \codesize
  \begin{center}
    \begin{code}
// Iteration over all vertices of all cells
unsigned int D = mesh.topology().dim();
for (MeshEntityIterator c(mesh, D); !c.end(); ++c)
{
  cout << "cell index = " << cell->index() << endl;
  for (MeshEntityIterator v(*c, 0); !v.end(); ++v)
  {
    cout << "vertex index = " << v->index() << endl;
  }
}

// Iteration over all vertices of all cells
for (CellIterator c(mesh); !c.end(); ++c)
{
  cout << "cell index = " << c->index() << endl;
  for (VertexIterator v(*c); !v.end(); ++v)
  {
    cout << "vertex index = " << v->index() << endl;
    cout << "vertex coordinates = " << v->point() << endl;
  }
}
    \end{code}
    \caption{Iteration over all vertices of all cells in a mesh, using
      the general iterator interface \texttt{MeshEntityIterator} and
      the specialized iterators \texttt{CellIterator} and \texttt{VertexIterator}.}
    \label{fig:interface:iterators}
  \end{center}
\end{figure}

\subsection{Direct Access to Mesh Data}

In addition to the iterator interface, all mesh data may be
accessed directly. Thus, one may obtain an array of the vertices of
all cells in the mesh directly from the mesh topology, and one may
obtain the vertex coordinates of the mesh directly from the mesh
geometry. This illustrated in Figure~\ref{fig:interface:direct} where
the same iteration as in Figure~\ref{fig:interface:iterators} is
performed without mesh iterators.

\begin{figure}[htbp]
  \codesize
  \begin{center}
    \begin{code}
MeshTopology& topology = mesh.topology();
MeshGeometry& geometry = mesh.geometry();
unsigned int D = topology.dim();
MeshConnectivity& connectivity = topology(D, 0);

for (unsigned int c = 0; c < topology.size(D); ++c)
{
  cout << "cell index = " << c << endl;

  unsigned int* vertices = connectivity(c);
  for (unsigned int i = 0; i < connectivity.size(c); ++i)
  {
    unsigned int vertex = vertices[i];

    cout << "vertex index = " << vertex << endl;
    cout << "vertex coordinates = "
         << geometry.point(vertex) << endl;
  }
}
    \end{code}
    \caption{Iteration over all vertices of all cells in a mesh and
      direct access of mesh data corresponding to the iteration of
      Figure~\ref{fig:interface:iterators} and Figure~\ref{fig:iteration}.}
    \label{fig:interface:direct}
  \end{center}
\end{figure}

\subsection{Mesh Algorithms}

In addition to the computation of mesh connectivity as discussed
previously, the DOLFIN mesh library provides a number of other useful
mesh algorithms, including boundary extraction, uniform mesh
refinement, adaptive mesh refinement (in preparation), mesh smoothing,
and reordering of mesh entities.

Figure~\ref{fig:interface:boundary} demonstrates uniform refinement
and boundary mesh extraction. When extracting a boundary mesh, it may
be desirable to also generate a mapping from the entities of the
boundary mesh to the corresponding entities of the original mesh. This
is the case for example when assembling the contribution from boundary
integrals during assembly of a linear system arising from a finite
element variational formulation of a PDE. One then needs to map each
cell of the boundary mesh to the corresponding facet of the original
mesh. (Note that the cells of the boundary mesh are facets of the
original mesh.) In Figure~\ref{fig:interface:boundary,meshfunctions},
we demonstrate how to extract a boundary and generate the mapping from
the boundary mesh to the original mesh. The mapping is expressed as
two \texttt{MeshFunction}s, one from the vertices of the boundary mesh
to the corresponding vertex indices of the original mesh and one from
the cells of the boundary mesh to the corresponding facet indices of
the original mesh.

\begin{figure}[htbp]
  \codesize
  \begin{center}
    \begin{code}
// Refine mesh uniformly twice
mesh.refine();
mesh.refine();

// Extract boundary mesh
BoundaryMesh boundary(mesh);

// Refine boundary mesh uniformly
boundary.refine();

// Save boundary mesh to file
File file("boundary.xml");
file << boundary;
    \end{code}
    \caption{Uniform refinement, boundary extraction and uniform
      refinement of the boundary mesh using the DOLFIN mesh
      library. Note that the extracted boundary mesh is itself a mesh
      and may thus for example be refined.}
    \label{fig:interface:boundary}
  \end{center}
\end{figure}

\begin{figure}[htbp]
  \codesize
  \begin{center}
    \begin{code}
MeshFunction<unsigned int> vertex_map;
MeshFunction<unsigned int> cell_map;

BoundaryMesh boundary(mesh, vertex_map, cell_map);
    \end{code}
    \caption{Extraction of a boundary mesh and generation of a pair of
    mappings from the vertices of the boundary mesh to the indices of
    the corresponding vertices of the original mesh and from the cells
    of the boundary mesh to the indices of the corresponding facets of
    the original mesh.}
    \label{fig:interface:boundary,meshfunctions}
  \end{center}
\end{figure}

\section{Parallel Considerations}

We discuss here how the concepts and data structures discussed above
can be used to assemble a global sparse finite element operator
(typically matrix) over a mesh distributed over several processors. It
is demonstrated below that we may reuse the concepts introduced above
to distribute the mesh. In particular, each processor owns a separate
piece of the global mesh, which can be stored as a regular
\texttt{Mesh}. Furthermore, each processor knows which facets of the
local mesh are incident with which facets on other processors and this
information can be stored as a pair of \texttt{MeshFunction}s.

\subsection{Simple Distribution of Mesh Data}

Let a mesh~$\mathcal{T}$ be given and assume that the mesh has been
partitioned into $n$ disjoint meshes $\{\mathcal{T}\}_{i=0}^{n-1}$
that together cover the computational domain $\Omega \subset
\R^n$. Such a partition can be computed using for example SCOTCH,
see~\cite{Pel2004}, or Metis, see~\cite{KarVip98a,KarVip98b}. On
each processor $p_i$, $i = 0,2,\ldots,n-1$, we store its part of the
global mesh as a regular mesh and in addition two mesh functions
$\mathcal{S}_i$ and $\mathcal{F}_i$ over the facets
of~$\mathcal{T}_i$.

We thus propose to store a distributed mesh~$\mathcal{T}$ on a set of
processors as the set of tuples $\{(\mathcal{T}_i, \mathcal{S}_i,
\mathcal{F}_i)\}_{i=0}^{n-1}$ where one tuple $(\mathcal{T}_i,
\mathcal{S}_i, \mathcal{F}_i)$ is stored on each processor $p_i$ for
$i=0,1,\ldots,n-1$.

The mesh function $\mathcal{S}_i$ maps each facet~$f$ to an integer $j
= \mathcal{S}_i(f)$ which indicates which (other) subdomain/mesh
$\mathcal{T}_j$ that the facet~$f$ is (physically) incident with,
\begin{equation}
  \mathcal{S}_i : (D-1, [0,1,\ldots,N^i_{D-1}-1]) \rightarrow [0,1,\ldots,n-1].
\end{equation}
Here, $(D-1, [0,1,\ldots,N^i_{D-1}-1])$ indicates that the domain of
$\mathcal{S}_i$ is the set of tuples $\{(D-1, k)\}$ where $0\leq k \leq
N^i_{D-1}-1$ and $N^i_{D-1}$ is the number of facets of
$\mathcal{T}_i$. Thus, if $j = \mathcal{S}_i(f)$ for some $j \neq i$,
then the facet~$f$ is shared with the mesh $\mathcal{T}_j$. If the
facet~$f$ is only incident with $\mathcal{T}_i$ itself, then we
set~$\mathcal{S}_i(f) = i$.

\begin{figure}[htbp]
  \begin{center}
    \includegraphics[width=8.5cm]{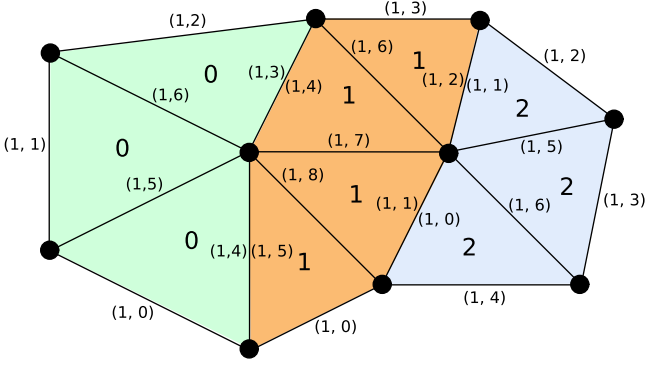}
    \caption{A global mesh~$\mathcal{T}$ partitioned into $n = 3$
      meshes $\{\mathcal{T}_i\}_{i=0}^{n-1}$. The mesh functions
      $\mathcal{S}_i$ and $\mathcal{F}_i$ indicate which facets of
      $\mathcal{T}_i$ are shared with other meshes. In this example,
      we have $\mathcal{S}_0((1, 3)) = \mathcal{S}_0((1, 4)) = 1$,
      indicating that facets $(1, 3)$ and $(1, 4)$ in $\mathcal{T}_0$
      are shared with $\mathcal{T}_1$. Furthermore, we may evaluate
      $\mathcal{F}_0$ at these two facets to find that the facet $(1,
      3)$ in $\mathcal{T}_0$ is incident to facet $(1, 4) =
      \mathcal{F}_0((1, 3))$ in $\mathcal{T}_1$ and facet $(1, 4)$ in
      $\mathcal{T}_0$ is incident to facet $(1, 5) = \mathcal{F}_0((1,
      4))$ in $\mathcal{T}_1$.}
    \label{fig:parallel}
  \end{center}
\end{figure}

The mesh function $\mathcal{F}_i$ maps each facet entity~$f$ to an
integer $\mathcal{F}_i(f)$ which indicates which facet~$f' = (D-1,
\mathcal{F}_i(f))$ of $\mathcal{T}_j$ for $j = \mathcal{S}_i(f)$ that
the facet~$f$ is incident (identical) to. If $\mathcal{S}_i(f) = i$,
then~$f$ is not shared with another mesh and we set $\mathcal{F}_i(f) = 0$.
 We illustrate the meaning of the two mesh
functions~$\mathcal{S}_i$ and $\mathcal{F}_i$ in
Figure~\ref{fig:parallel}.

\subsection{Parallel Assembly}

The standard algorithm for computing a global sparse operator (tensor)
from a finite element variational form is known as \emph{assembly},
see~\cite{ZieTay67,Hug87,Lan99}. By this algorithm, the global sparse
operator may be computed by assembling (summing) the contributions
from the local entities of a finite element mesh. On each cell of the
mesh, one computes a small \emph{cell tensor} (often referred to as
the ``element stiffness matrix'') and add the entries of that tensor
to a global sparse tensor (often referred to as the ``global stiffness
matrix''). We shall not discuss the assembly algorithm in detail here
and refer instead to~\cite{AlnLan2008,logg:article:12,logg:manual:04}, but note
that to add the entries from the local cell tensor to the global
sparse tensor, we need to compute a so-called \emph{local-to-global}
mapping on each cell. This maps the local degrees of freedom on a cell
(numbering the rows and columns of the cell tensor) to global degrees
of freedom (numbering the rows and columns of the global tensor).

The assembly algorithm may be trivially parallelized by letting each
processor $p_i$ compute and insert the cell tensors on the local mesh
$\mathcal{T}_i$ into the global tensor. When the global tensor is a
sparse matrix, linear algebra libraries like PETSc,
see~\cite{BalBus04,BalGro97}, may be used to store the global tensor
in parallel. PETSc handles the communication of matrix data between
processors and we need only make sure that each processor knows how to
insert entries into the global sparse matrix according to the
local-to-global mapping of the global mesh. We demonstrate below that
on each processor $p_i$, we may (with a small amount of communication
with neighboring processors) compute the part of the local-to-global
mapping of the global mesh~$\mathcal{T}$ relevant to each local
mesh~$\mathcal{T}_i$ in parallel on each processor~$p_i$, which thus
allows us to assemble the global sparse matrix in parallel.

\subsection{Mapping Degrees of Freedom in Parallel}

In Algorithm~\ref{alg:dofmap}, we describe how the mapping of degrees
of freedom may be computed on a distributed mesh $\{(\mathcal{T}_i,
\mathcal{S}_i, \mathcal{F}_i)\}_{i=0}^{n-1}$.

To express this algorithm in compact form, we need to introduce some
further notation. For each cell~$c$ in a local finite element
mesh~$\mathcal{T}_i$, we assume that we can compute a local-to-global
mapping $\iota^i_c$, which maps each local degree of freedom on the
cell~$c$ to a global degree of freedom (for a numbering scheme valid
on the local mesh~$\mathcal{T}_i$). For example, when computing with
standard piecewise linear finite elements on triangles, the
local-to-global mapping~$\iota^i_c$ may map a local vertex number~$0$,
$1$ or $2$ on~$c$ to the corresponding global number vertex number on
the mesh~$\mathcal{T}_i$. Thus, the domain of $\iota^i_c$ is here
$\{0,1,2\}$ and the range is $[0,N^i_0-1]$, where $N^i_0$ is the
number of vertices of the mesh~$\mathcal{T}_i$. We emphasize that the
local-to-global mapping~$\iota^i_c$ is not aware of the global
mesh~$\mathcal{T}$ of which the local mesh $\mathcal{T}_i$ is a part.
Instead, it is the task of Algorithm~\ref{alg:dofmap} to compute (in
parallel) a local-to-global mapping valid on the global
mesh~$\mathcal{T}$ from a given local-to-global mapping on each local
mesh~$\mathcal{T}_i$.

We let~$\mathcal{M}_i$ denote the parallel local-to-global mapping to
be computed on each part~$\mathcal{T}_i$ of the global mesh. For ease
of notation, we express $\mathcal{M}_i$ as a set of tuples
$\mathcal{M}_i = \{((c,i), I)\}$, where~$c$ is a cell (number), $i$ is
the local number of a degree of freedom on~$c$ and $I$ is the
corresponding global number. The mapping~$\mathcal{M}_i$ should be
thought of as a function that maps a cell and a local degree of
freedom~$(c, i)$ to the corresponding global degree of freedom~$I$. In
the case of standard piecewise linear elements on triangles, the
domain of $\mathcal{M}_i$ is $[0,N^i_D-1] \times \{0,1,2\}$, where
$N^i_D$ is the number of cells of the mesh~$\mathcal{T}_i$, and the
range of $\mathcal{M}_i$ is $[0,N_0-1]$, where $N_0$ is the total
number of vertices of the partitioned global mesh~$\mathcal{T}$. We
note here that the mapping $\mathcal{M}_i$ may be stored as a
fixed-size array.

To compute the parallel local-to-global mapping~$\mathcal{M}_i$ on
each processor, we need to iterate over the entities of the
mesh~$\mathcal{T}_i$ and renumber the degrees of freedom. To do this,
we introduce an auxiliary temporary mapping~$\mathcal{N}_i$ on each
processor that maps degrees of freedom as given by the local-to-global
mapping~$\iota^i_c$ on each cell~$c$ (which is only aware of how to map
degrees of freedom internally on~$\mathcal{T}_i$) to degrees of
freedom as given by the parallel local-to-global
mapping~$\mathcal{M}_i$ (which is aware of how to map degrees of
freedom globally on the distributed mesh). In the case of standard
piecewise linear elements on triangles, the domain of~$\mathcal{N}_i$
is $[0,N^i_0-1]$ and the range of~$\mathcal{N}_i$ is $[0,N_0-1]$.
Just as for $\mathcal{M}_i$, we may think of $\mathcal{N}_i$ as a
function but write it as a set of tuples (pairs)
in~Algorithm~\ref{alg:dofmap} for ease of notation. In a C++
implementation, $\mathcal{N}_i$ may be stored in the form of an
STL map (\texttt{std::map<unsigned int, unsigned int>}).

Finally, we let $\overline{N}_i$ denote the number of degrees of
freedom on~$\mathcal{T}_i$ not shared with a mesh $\mathcal{T}_j$ for
$j < i$. This number can be computed on each mesh~$\mathcal{T}_i$ from
the mesh function~$\mathcal{S}_i$.

In Algorithm~\ref{alg:dofmap}, we first let each processor $p_i$
compute $\overline{N}_i$ (in parallel). These numbers are then
communicated successively from $p_{i-1}$ to $p_i$ to compute an offset
for the numbering of degrees of freedom on each mesh $\mathcal{T}_i$.
We then let each processor $p_i$ number the degrees of freedom (in
parallel) on cells which are incident with the boundary of
$\mathcal{T}_i$ and are \emph{shared with a mesh~$\mathcal{T}_j$
for $j > i$}. After the degrees of freedom on mesh boundaries have
been numbered on each processor, those numbering schemes are
communicated successively from $p_i$ to $p_j$ for all $i<j$
such that $\mathcal{T}_i$ and $\mathcal{T}_j$ share degrees of freedom
on a common facet. Finally, each processor~$p_i$ may number the
$\overline{N}_i$ ``internal degrees of freedom'' on $\mathcal{T}_i$.

The key point of this algorithm is to always let $p_i$ number degrees
of freedom common with $p_j$ for $i < j$. This numbering of shared
degrees of freedom is then communicated over shared facets, from $f'$
to $f$ in stage~2 of Algorithm~\ref{alg:dofmap}. Since it is known
a~priori which facets are shared between two meshes $\mathcal{T}_i$ and
$\mathcal{T}_j$, one may communicate the common numbering for all
shared degrees of freedom from $p_i$ to $p_j$ in one batch.

\begin{figure}[htbp]
  \begin{center}
    \includegraphics[width=8.5cm]{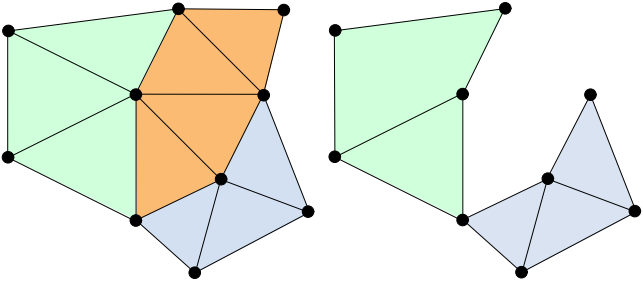}
    \caption{Two partitions of a mesh. In the partition on the left,
    two of the meshes share only a common vertex and may thus share a
    single degree of freedom at that vertex. The communication of a
    common numbering of that degree of freedom is propagated over the
    facets incident with the common vertex to the neighboring mesh.
    In the partition on the right, it is not possible to propagate the
    numbering over facets (since there are no shared facets) and so
    Algorithm~\ref{alg:dofmap} will fail to compute a correct
    numbering scheme for this partition.}
    \label{fig:allowed}
  \end{center}
\end{figure}

Note that it is important that the communication of shared degrees of
freedom in stage~2 of Algorithm~\ref{alg:dofmap} is performed
sequentially, starting with processor $p_1$ receiving the common
numbering from $p_0$, then $p_2$ receiving the common numbering from
$p_0$ and/or $p_1$ etc. This guarantees that common degrees of freedom
are communicated from $p_i$ to all $p_j$ with $i < j$ such that
$\mathcal{T}_i$ and $\mathcal{T}_j$ share common degrees of freedom,
even if $\mathcal{T}_i$ and $\mathcal{T}_j$ don't share a common
facet. For this to work, we make the assumption that if any two meshes
$\mathcal{T}_i$ and $\mathcal{T}_j$ share a common degree of freedom
on the two cells~$c\in\mathcal{T}_i$ and $c'\in\mathcal{T}_j$, then
each of $\mathcal{T}_i$ and $\mathcal{T}_j$ must share that degree of
freedom with some other mesh $\mathcal{T}_{i'}$ or $\mathcal{T}_{j'}$
respectively. In Figure~\ref{fig:allowed}, we illustrate this
assumption (for $i' = j'$) and give an example of a partition for
which Algorithm~\ref{alg:dofmap} fails to correctly number all shared
degrees of freedom. It is a mild assumption to disallow such
partitions (and meshes).

Algorithm~\ref{alg:dofmap} is currently not implemented in
DOLFIN. Instead, a simple (but suboptimal) strategy where the
computational mesh is broadcast to all processors has been
implemented. Each processor owns a copy of the entire mesh and knows
which part of that mesh to assemble. For a further discussion of the
current implementation of parallel assembly in DOLFIN (available with
DOLFIN 0.7.2), see~\cite{Vik2008}.

\begin{algorithm}[htbp]
  \begin{tabbing}
    --- Stage 0: \emph{Compute offsets} \\
    \textbf{on each} processor $p_i$ \\
    \tab compute $\overline{N}_i$ \\
    \textbf{on} processor $p_0$ \\
    \tab $\mathrm{offset}_0 = 0$ \\
    \textbf{for} $i = 1,2,\ldots,n-1$ \\
    \tab \textbf{on} processor $p_i$ \\
    \tab \tab Receive $(\mathrm{offset}_{i-1}, \overline{N}_{i-1})$ from $p_{i-1}$ \\
    \tab \tab $\mathrm{offset}_i = \mathrm{offset}_{i-1} + \overline{N}_{i-1}$ \\
    --- Stage 1: \emph{Compute mapping on shared facets} \\
    \textbf{on each} processor $p_i$ \\
    \tab $\mathcal{M}_i = \emptyset$ \\
    \tab $\mathcal{N}_i = \emptyset$ \\
    \tab $k_i = \overline{N}_i$ \\
    \tab \textbf{for each} facet $f \in \mathcal{T}_i$ \\
    \tab \tab $j = \mathcal{S}_i(f)$ \\
    \tab \tab \textbf{if} $j > i$ \\
    \tab \tab \tab Let $c\in\mathcal{T}_i$ be the cell incident with $f$ \\
    \tab \tab \tab \textbf{for each} local degree of freedom $l$ on $c$ \\
    \tab \tab \tab \tab \textbf{if} $(\iota^i_c(l), L) \in \mathcal{N}_i$ for some $L$ \\
    \tab \tab \tab \tab \tab $\mathcal{M}_i = \mathcal{M}_i \cup ((c, l), L)$ \\
    \tab \tab \tab \tab \textbf{else} \\
    \tab \tab \tab \tab \tab $\mathcal{M}_i = \mathcal{M}_i \cup ((c, l), k_i)$ \\
    \tab \tab \tab \tab \tab $\mathcal{N}_i = \mathcal{N}_i \cup (\iota^i_c(l), k_i)$ \\
    \tab \tab \tab \tab \tab $k_i = k_i + 1$ \\
    --- Stage 2: \emph{Communicate mapping on shared facets} \\
    \textbf{for} $j=1,2,\ldots,n-1$ \\
    \tab \textbf{for each} facet $f \in \mathcal{T}_j$ \\
    \tab \tab $i = \mathcal{S}_j(f)$ \\
    \tab \tab \textbf{if} $i < j$ \\
    \tab \tab \tab $f' = \mathcal{F}_j(f)$ \\
    \tab \tab \tab Receive degrees of freedom on $f'$ from $p_i$ \\
    \tab \tab \tab Update $\mathcal{N}_j$ for shared degrees of freedom \\
    --- Stage 3: \emph{Compute mapping for interior degrees of freedom} \\
    \textbf{on each} processor $p_i$ \\
    \tab \textbf{for each} cell $c\in\mathcal{T}_i$ \\
    \tab \tab \textbf{for each} local degree of freedom $l$ on $c$ \\
    \tab \tab \tab \textbf{if} $(\iota^i_c(l), L) \in \mathcal{N}_i$ for some $L$ \\
    \tab \tab \tab \tab $\mathcal{M}_i = \mathcal{M}_i \cup ((c, l), L)$ \\
    \tab \tab \tab \textbf{else} \\
    \tab \tab \tab \tab $\mathcal{M}_i = \mathcal{M}_i \cup ((c, l), k_i)$ \\
    \tab \tab \tab \tab $\mathcal{N}_i = \mathcal{N}_i \cup (\iota^i_c(l), k_i)$ \\
    \tab \tab \tab \tab $k_i = k_i + 1$ \\
  \end{tabbing}
  \caption{$\{\mathcal{M}_i\}$ = ComputeMapping($\{(\mathcal{T}_i, \mathcal{S}_i,
    \mathcal{F}_i)\}$), computing the local-to-global mapping
    in parallel for a mesh distributed over $n$ processors $p_i$,
    $i=0,1,\ldots,n-1$.}
  \label{alg:dofmap}
\end{algorithm}

\section{Benchmark Results}

In this section, we present a series of benchmarks to illustrate the
efficiency of the mesh representation and its implementation. The new
mesh library (which is available as part of DOLFIN since version
0.6.3) is compared to the old DOLFIN mesh library which is a fairly
efficient C++ implementation, but which suffers from object-oriented
overhead; all mesh entities are there stored as arrays of objects,
which store their data locally in each object (including mesh
incidence relations and vertex coordinates). The benchmark results
were obtained on a 2.66~GHz 64-bit Intel processor (Q6700) running
Ubuntu GNU/Linux for DOLFIN version 0.6.2-1 and DOLFIN version 0.7.1
respectively.

The five test cases that are examined are the following: (i) CPU time
and (ii) memory usage for creation of a uniform tetrahedral mesh of
the unit cube, (iii) CPU time for iteration over all vertices of the
mesh, (iv) CPU time for accessing the coordinates of all vertices of
the mesh, and (iv) uniform refinement of the mesh.

In summary, the speedup was in all cases a factor $10$--$1000$ and
memory usage was reduced by a factor of~$50$ (comparing slopes of
lines in Figure~\ref{fig:bench,memory}). The speedup and decreased
memory usage is the result of more efficient algorithms and data
structures, where all data is stored in large static arrays and
objects are only provided as part of the interface for simple access
to the underlying data representation, not to store data
themselves. Another contributing factor is that the old DOLFIN mesh
library precomputes certain connectivities (including the edges and
faces of each cell) at startup, whereas this computation is carried
out only when requested in the new DOLFIN mesh library, either as part
of the iterator interface or by an explicit call to
\texttt{Mesh::init()}. We also note that from
Figure~\ref{fig:bench,memory}, it is evident that DOLFIN can be
further improved in terms of its minimal memory footprint.

\begin{figure}[htbp]
  \begin{center}
    \includegraphics[width=9cm]{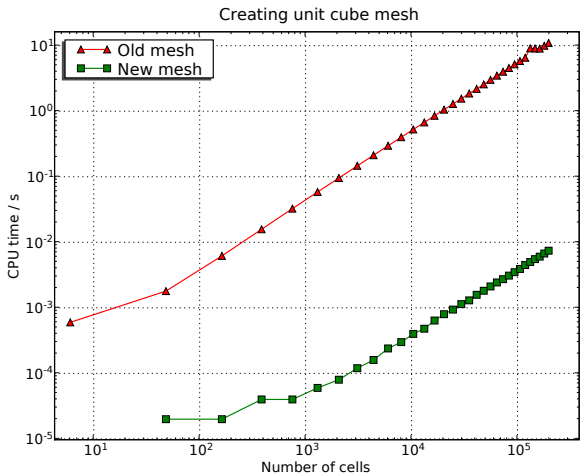}
    \caption{Benchmarking the CPU time for creation of a uniform tetrahedral mesh of the unit cube
             for the new mesh library vs. the old DOLFIN mesh library.}
    \label{fig:bench,create}
  \end{center}
\end{figure}

\begin{figure}[htbp]
  \begin{center}
  \includegraphics[width=9cm]{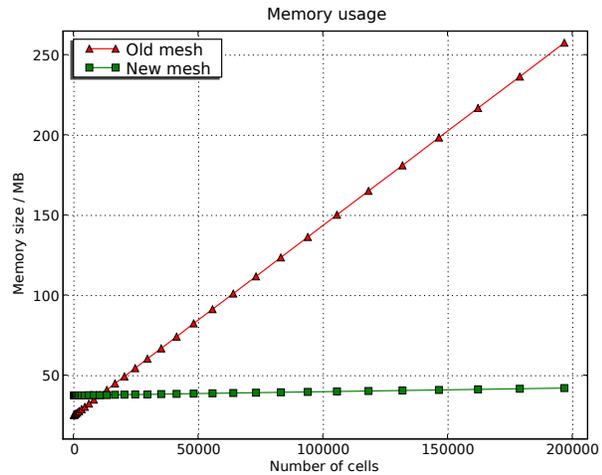}
    \caption{Benchmarking the memory usage for creation of a uniform tetrahedral mesh of the unit cube
             for the new mesh library vs. the old DOLFIN mesh library.}
    \label{fig:bench,memory}
  \end{center}
\end{figure}

\begin{figure}[htbp]
  \begin{center}
    \includegraphics[width=9cm]{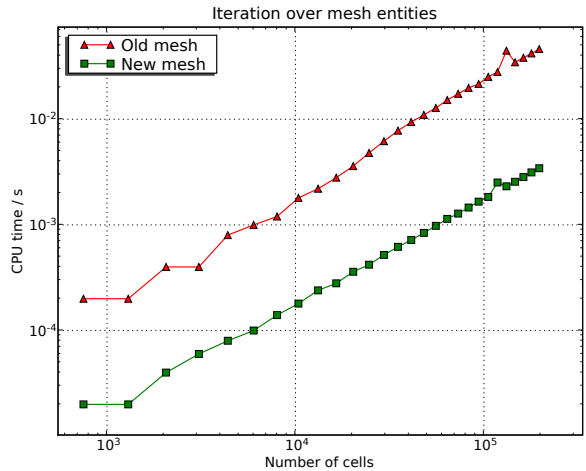}
    \caption{Benchmarking the CPU time for iteration over the vertices of each cell
             for the new mesh library vs. the old DOLFIN mesh library.}
    \label{fig:bench,entities}
  \end{center}
\end{figure}

\begin{figure}[htbp]
  \begin{center}
    \includegraphics[width=9cm]{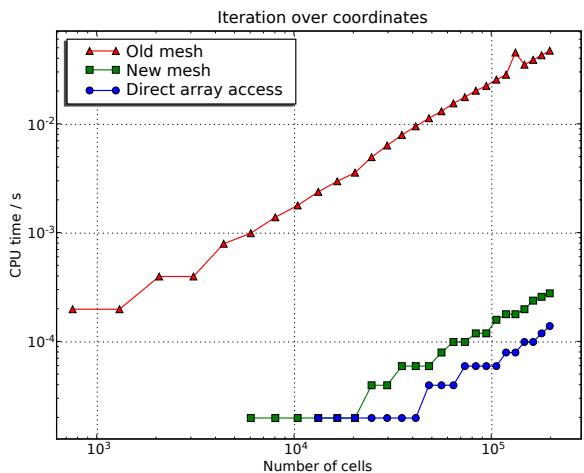}
    \caption{Benchmarking the CPU time for iteration over the coordinates of each vertex
             for the new mesh library vs. the old mesh library vs.
             direct access of coordinate arrays in the new DOLFIN mesh library.}
    \label{fig:bench,coordinates}
  \end{center}
\end{figure}

\begin{figure}[htbp]
  \begin{center}
  \includegraphics[width=9cm]{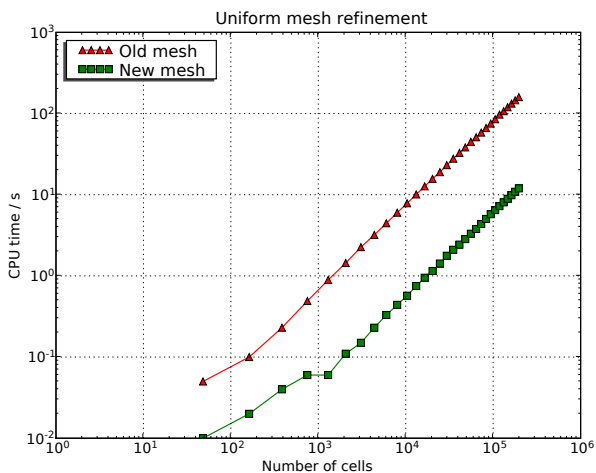}
    \caption{Benchmarking the CPU time for uniform refinement of the unit cube
             for the new mesh library vs. the old DOLFIN mesh library.}
    \label{fig:bench,refinement}
  \end{center}
\end{figure}

\section{Conclusions}

We have presented a simple yet general and efficient representation of
computational meshes and demonstrated a straightforward implementation
of this representation as a set of C++ classes that correspond to the
basic concepts of the mesh representation. The implementation is
available freely as part of DOLFIN~\cite{www:dolfin}.

\section*{Acknowledgments}

The author wishes to acknowledge the many contributions to the DOLFIN
mesh library from the DOLFIN developers, in particular Garth N. Wells,
Johan Hoffman, Johan Jansson, Kristian Oelgaard, Ola Skavhaug and
Gustav Magnus Vikstr\o{}m. The author also wishes to acknowledge
invaluable input from and inspiring discussions with the authors
of~\cite{KnepleyKarpeev07A}.

\vspace{20mm}

\bibliographystyle{siam}
\bibliography{bibliography}

\end{document}